\documentclass [12pt,a4paper,reqno]{amsart}
\usepackage{amsmath,amsthm}
\usepackage{amscd}
\usepackage{amsfonts}
\input amssymb.sty
\chardef\bslash=`\\ 

\newtheorem{thm}{Theorem}
\newtheorem*{thm*}{Theorem}
\newtheorem*{conjecture*}{Conjecture}

\newtheorem{cor}{Corollary}

\newtheorem{lem}{Lemma}

\theoremstyle{definition}
\newtheorem{defn}{Definition}
\newtheorem*{remark*}{Remarks}
\newtheorem*{examples*}{Examples}
\newtheorem*{defn*}{Definition}
\newtheorem*{cor*}{Corollary}

\newcommand{\thmref}[1]{Theorem~\ref{#1}}
\newcommand{\secref}[1]{Section~\ref{#1}}

\newcommand{\lemref}[1]{Lemma~\ref{#1}}
\newcommand{\corref}[1]{Corollary~\ref{#1}}

\numberwithin{equation}{section}

\newcommand{\Th}{\Theta}
\newcommand{\C}{\mathcal{C}}
\newcommand{\D}{\mathcal{D}}
\newcommand{\E}{\mathcal{E}}
\newcommand{\V}{\mathcal{V}}
\newcommand{\Q}{\mathcal{Q}}
\newcommand{\st}{\sigma}

 \renewcommand{\sectionmark}[1]{}

\newcommand{\field}[1]{\mathbb{#1}}

\newcommand{\vp}{\varphi}

\newcommand{\Z}{\field{Z}}
\newcommand{\K}{\mathbf{K}}

\renewcommand{\a}{\alpha}

\date{}
\begin{document}

\title{On automorphisms of categories of free algebras of some varieties}

\author[Boris Plotkin]{B. Plotkin\\
 Institute of Mathematics \\
 Hebrew University, 91803 Jerusalem, Israel\\}

\author[Grigori Zhitomirski]{G. Zhitomirski\\
 Department of Mathematics\\
 Bar-Ilan University, 52900 Ramat Gan, Israel}
\thanks {This research of the second author was partially supported by THE ISRAEL SCIENCE FOUNDATION
 founded by The Israel Academy of Sciences and Humanities - Center of
 Excellence Program.}

\maketitle

\begin{abstract}
Let  $\V$ be a variety of universal algebras. A method is suggested for describing automorphisms of a category
of free $\V$-algebras. Applying this general method all automorphisms of such categories are found in two cases:
1) $\V$ is the variety of all free associative $K-$algebras over an infinite field $\K$ and 2) $\V$ is the
variety of all representations of groups in unital $R-$modules over a commutative associative ring with unit. It
turns out that they are almost inner in a sense.
\end{abstract}

\baselineskip 20pt
\bigskip

\section{INTRODUCTION}\label{Intro}

Let $\V$ be a variety of universal algebras. Consider the category $\Th (\V)$ whose objects are all algebras
from $\V$ and whose morphisms are all homomorphisms of them. Fix an infinite set $X\sb 0$. Let $\Th \sp 0 (\V )$
be the full subcategory of $\Th (\V )$ defined by all free algebras from $\V$ over finite subsets of the set
$X\sb 0$. The problem is to describe all automorphisms of the category $\Th \sp 0 (\V )$. Motivations for this
research direction can be found in papers \cite {Free}, \cite{LieAlg}, \cite{AlgGeom}. We only remind here that
it is connected with some problems in universal algebraic geometry, particularly, when two algebras of the
considered variety have the same geometry.

The mentioned problem was solved in many cases, for example, for varieties of all groups and all semigroups
\cite{Free}, inverse semigroups \cite{InvSem} and for the variety of all Lie algebras \cite{LieAlg}. In all
known cases the automorphisms turn out to be inner or close to inner (almost inner).

Recall that an automorphism $\Phi$ of a category $\C$ is called inner if it is isomorphic to the identity
functor $Id \sp {\C}$   in the category of all endofunctors of $\C$. It means that there exists a function that
assigns to every object $A$ of the given category an isomorphism $\st \sb A : A\to \Phi (A)$ such that for every
morphism $\mu: A\to B$ we have $\Phi (\mu)= \st \sb B \circ \mu \circ \st \sb A \sp {-1}$. This fact explains
the term "inner". Thus if an automorphism $\Phi$ is inner the object $\Phi (A)$ is isomorphic to $A$ for every
$\C$-object $A$.

For most interesting varieties, the last condition is satisfied. Thus there is to assume that every
automorphisms $\Phi$ of the category $\Th \sp 0 (\V )$ takes every its object to an isomorphic one. It was
proved \cite {Free} that in this case $\Phi$ is a product of two automorphisms one of each is inner and the
other one is stable, that is, it leaves fixed all objects. Every stable automorphism, for its part, induces an
automorphism of the monoid $End A$ for every object $A$. Therefore the usual method to describe automorphisms of
the category $\Th \sp 0 (\V )$ includes describing of $AutEnd A$ for free finitely generated algebras in the
variety $\V$. But it turns out that the last problem, being interesting itself, is rather complicated than the
original one (see \cite{Form} ,\cite{SemFree}, \cite{InvSem}).

Developing ideas and results of the paper \cite{PlZh}, we suggest a new method of describing automorphisms of
categories of universal algebras. Below, we explain this method.

1. Let $A\sb 0$ be a monogenic free algebra in a category $\C$ of universal algebras and let $x\sb 0$ be its
free generator. Let $A$ be a $\C-$algebra and $a\in A$. Denote by $\a \sb a$ the unique homomorphism from $A\sb
0$ to $A$ that takes $x\sb 0$ to $a$: $\a \sb a (x\sb 0) = a$. It is clear that it is an one-to-one
correspondence between sets $A$ and $Hom(A\sb 0,A)$: $a\mapsto \a \sb a$. Thus every automorphism $\Phi$ of the
category $\C$ leaving fixed $A\sb 0$ determines a family of bijections $(s\sp {\Phi} \sb A\;\vert A\in Ob\,\C )$
defined by the following way: $s\sp {\Phi} \sb A (a)=\Phi (\a \sb a )(x\sb 0)$ for every $a\in A$.

It is easy to see that if $\Phi$ is an identity automorphism then $s\sp {\Phi} \sb A =1\sb A$,  if $\Phi$ is a
product $\Gamma \circ \Psi$ then $s\sp {\Phi} \sb A =s\sp {\Gamma} \sb A \circ s\sp {\Psi} \sb A $ and if $\Psi
=\Phi \sp {-1}$ then $s\sp {\Psi} \sb A =(s\sp {\Phi} \sb A)\sp {-1}$. It follows from the definition of the
function $A\mapsto s\sp {\Phi} \sb A $  that for every homomorphism $\nu :A\to B$, where $A,B$ are objects of
$\C$, we have $\Phi (\nu)=s \sb B \circ \nu \circ s \sb A \sp {-1}$. The last fact leads to idea to introduce a
notion of potential-inner automorphism.

2. Let $\D$ be an extension of a category $\C$ obtained by adding some new maps as morphisms. An automorphism
$\Phi$ of the category $\C$ is called $\D-inner$ if it is the restriction to the category $\C$ of some inner
automorphism of the category, that is, there exists a function $f$ which assigns to every object $A$ a
$\D-$isomorphism $f\sb A : A\to \Phi (A)$ such that $\Phi (\mu)= f\sb B \circ \mu \circ f\sb A \sp {-1}$ holds
for every $\C$-morphism $\mu: A\to B$. An automorphism $\Phi $ of a given category $\C$ is called potential
inner if it is $\D$-inner for some extension $\D$ of $\C$.

This notion gives an opportunity to consider the problem from a new point of view, namely it can be formulated
now in the following way: 1) what extension $\D$ of the given category $\C$ we have to construct in order to
make all $\C$-automorphisms to be $\D$-inner and 2) when potential-inner automorphisms are in fact inner.

It turns out that an automorphism $\Phi $ of a category $\C$ containing a monogenic free algebra $A\sb 0$ is
potential-inner if and only if $\Phi (A\sb 0)$ is isomorphic to $A\sb 0$. (\thmref{potent-inner}). It was
mentioned above that we may consider only the case of potential-inner automorphisms.

3. It is shown (\lemref{subclass}) that we can assume that $\Phi $ leaves fixed all free algebras in $\C$. Using
bijections $s\sb A$ we can define a new algebraic structure $A\sp *$ on the underlying set of every free algebra
$A$  such that $s\sb A :A\to A\sp *$ is an isomorphism. These new algebras $A\sp *$ need not be the objects of
the category $\C$, but it turns out that $A\sp *$ can be described. It is the crucial point of our idea. If $A$
is a free algebra such that the number of its free generators is not less than arities of all its operations
then $A\sp *$ is a derivative algebraic structure, i.e. all its basic operations are determined by terms of
corresponding language or, in other words, are polynomial operations in algebra $A$ (\thmref{main}).

Thus every map $s\sb A$ is an isomorphism between the source structure and the derivative structure on the same
set. It is worth mentioning that the source structure is a derivative structure for $A\sp *$. It helps us to
find the derivative structures $A\sp *$ and to describe the maps $s\sb A$.

4. Describing maps $s\sb A$ leads to a description of a given automorphism, that is, we know its form. But it is
not a unique description, it may be that there exists a better one. It turns out (\lemref{InnerAndCentral}) that
the last problem can be solved by means of so called central functions.

In the next section we prove all facts mentioned above. Then we apply the suggested method and describe all
automorphisms of the category of free associative algebras (\secref{secasslinalg}) and all automorphisms of the
category of free group representations (\secref{secGroupRep}). In both cases the automorphisms turn out to be
almost inner. We also show that our method give an opportunity to obtain more simply many already known results.

Thus the main results of the paper are the following ones. \thmref{potent-inner} gives a necessary and
sufficient condition for an automorphism to be potential inner and very impotent property of such automorphisms.
\thmref{saveGenerators} reduces the problem to the case when maps $s\sb A$ leave fixed all basis elements of
every free algebra $A$. \thmref{main} describes maps $s\sb A$. \thmref{FreeAssAlg} characterizes automorphisms
of a category of free associative algebras as almost inner\footnote{Others proofs of this result are given by
Berzins \cite{Berz} and Mashevitzky \cite{MashPrep} using quite different approaches.}, and the similar result
is presented by (\thmref{Rep-R}) for categories of group representations.

All used and not defined notions of category theory and universal algebra can be found in \cite{Categ} and
\cite{Gratz}

\section{basic results}\label{basic}
\subsection{Objects}
In this section, we repeat some of general results from \cite{PlZh} changing them according to the case of
categories of free universal algebras only. We consider arbitrary universal algebras, maybe many-sorted, with
some notion of standard homomorphisms. In the case of one-sorted algebras the notion of a homomorphism is the
usual one and in the many-sorted case it may be defined in a special way, usually it is a finite sequence of
corresponding maps. We assume that we can speak about varieties of such algebras. We also assume that there is a
notion of free algebras in these varieties over arbitrary sets of free generators. Just in case, every map from
a set $X$ into a $\V -$algebra $A$ can be extended to a unique homomorphism of the free algebra over $X$ to $A$.

Given a variety $\V $. We denote by the $\Th (\V ) $ category whose objects are algebras from $\V $ and
morphisms are the standard homomorphisms of mentioned algebras. We assume also that we have a forgetful functor,
that is, a faithful functor from $\Th (\V ) $ to the category of all sets and maps. In an one-sorted case, this
functor assigns to every algebra its underlaying set and to every homomorphism itself as a map. A many-sorted
case is more difficult and demands a special approach.

\subsection{Inner and potential inner automorphisms}
All categories we consider are full subcategories of $\Th (\V ) $. We assume that all of them contain a
monogenic (one-generated) free algebra.  Because morphisms of $\Th (\V ) $ are maps we can consider an extension
$\Q $ of $\Th (\V ) $ with the same objects but extended sets of morphisms between them. The new morphisms added
in some way will be called quasi-homomorphisms. Let $\C $ be a full subcategory of the category $\Th (\V ) $.
{\it An extension } $\D $ of $\C $ will always mean a full subcategory of $\Q$ with the same objects that $\C $.
That is, we add to $\C $ some new morphisms, namely, quasi-homomorphisms from $\Q $.

\begin{defn}\label{D-inner}
An automorphism $\Phi$ of a category $\C $ is called {\it inner } if for every object $A$ of this category there
exists an isomorphism $\st \sb A : A\to \Phi (A)$  such that for every homomorphism $\mu: A\to B$ we have $\Phi
(\mu)= \st \sb B \circ \mu \circ \st \sb A \sp {-1}$. That is the following diagram is commutative:
$$
\CD
A @>\st \sb A>>    \Phi (A)\\
@V\mu VV          @VV \Phi (\mu) V\\
B @>\st \sb B  >>   \Phi (B).
\endCD
$$
In other words, $\Phi$ is isomorphic as a functor to the identity functor $Id$ of $\C$.

An automorphism $\Phi$  is called  $\D -${\it inner} if it is the restriction to $\C$ of some inner automorphism
of an extension $\D$ of the category $\C$ or, on other words, $\st -$s in the diagram above are $\Q
-$isomorphisms (quasi-isomorphisms).
\end{defn}

\begin{defn} An automorphism $\Phi$ of the category $\C $ is said to be  {\it potential-inner} if it is
$\D$-inner for some extension $\D$ of $\C $.
\end{defn}

Under hypothesis, there exists a free monogenic algebra in the considered categories.  Let $A\sb 0 $ be a free
monogenic algebra over a fixed element $x\sb 0$ in a category $\C$. There is a bijection between underlying set
of an $\C -$algebra $A$ and the set of all homomorphisms from $A\sb 0$ into $A$, namely, to every element $a\in
A$ corresponds the homomorphism $\a \sb a $ defined by:
\begin{equation}\label{alphas}
\a \sb a \sp A (x\sb 0 ) =a.
\end{equation}

The next result gives a necessary and sufficient condition for an automorphism of a category of universal
algebras to be potential inner.
\begin{thm}\label{potent-inner}  An automorphism $\Phi$ of the category $\C $ is potential-inner if and only if
$\Phi (A\sb 0)$ is isomorphic to $A\sb 0$.

\noindent If $\Phi $ is potential inner, i.e., there is a function $A \mapsto s\sb A$ such that $\Phi (\mu)= s
\sb B \circ \mu \circ s \sb A \sp {-1}$ for every  $\mu: A\to B$, then $\Phi (F)$ is isomorphic to $F$ for every
free algebra $F$ in the category $\C$ and quasi-isomorphism $s\sb F$ maps basis of $F$ onto basis of $\Phi (F)$.
\end{thm}
\begin{proof}
Let $\st :A\sb 0 \to \Phi (A\sb 0)$ be an isomorphism. Let $A$ be an arbitrary $\C -$algebra and $a\in A$. Using
the formula~\ref{alphas},  there exists one and only one element $\bar{a}\in \Phi (A)$ such that
 $\Phi (\a \sb a \sp A)\circ \st =\a \sb {\bar{a}}\sp {\Phi (A)}$. Since
 $\a \sb a \sp A =\Phi \sp {-1}(\a \sb {\bar{a}}\sp {\Phi (A)}\circ \st \sp {-1})$, we obtain a bijection
$s\sb A : A\to \Phi (A)$  setting for every $a\in A$:
\begin{equation}\label{s-maps}
s\sb A (a)= \Phi (\alpha \sb a \sp A )\circ \st (x\sb 0)
\end{equation}
and hence a family of bijections $(s\sb A : A\to \Phi (A),\; \vert A\in Ob\C )$.

Let $\nu :A\to B $ be a homomorphism. According to the definition above we have $s\sb B (\nu (a))= \Phi (\alpha
\sb {\nu (a)}\sp B)\circ \st (x\sb 0)$. Since $\alpha \sb {\nu (a)}\sp B=\nu \circ \alpha \sb a \sp A$, we
obtain that $(s\sb B \circ \nu) (a)=\Phi (\nu )\circ \Phi (\alpha \sb a \sp A)\circ \st (x\sb 0)=(\Phi (\nu
)\circ s \sb A )(a)$. Hence
\begin{equation}\label{auto-action}
 \Phi (\nu )=s \sb B \circ \nu \circ s \sb A \sp {-1}.
\end{equation}

Consider a category $\D$  containing $\C$ as a subcategory and generated by adding new isomorphisms (bijections)
$s\sb A: A\to \Phi (A)$ and their inverses $s\sp {-1}\sb A: \Phi (A)\to A $. Under definition the automorphism
$\Phi$ is $\D$-inner.

Let now an automorphism $\Phi $ of $\C$ is potential inner and $F$ be a free algebra in the category $\C$ over a
set $X$. Denote $A= \Phi (F)$ and $B= \Phi \sp {-1}(F)$. We have the following diagram:
$$
\CD B @> s\sb B>>    F @> s\sb F >> A
\endCD
$$
Set $\tilde{X}=s\sb F (X)$ and note that $s\sb F$ determines a bijection between $X$ and $\tilde{X}$. Denote by
$\st$ the unique homomorphism from $F$ into $A= \Phi (F)$ that extends $s\sb F \vert \sb X$, that is, $\st
(x)=s\sb F (x)$ for every $x\in X$. In the same way, we have homomorphism $\tau :F \to B$ such that $\tau
(x)=s\sp {-1}\sb B (x)$ for every $x\in X$ . Let us calculate:
$$\Phi (\tau) \circ \st (x)= s \sb B \circ \tau \circ s \sb F \sp {-1}\circ s\sb F (x)
=s \sb B \circ s\sb B \sp {-1} (x)=x .$$ Hence $\Phi (\tau) \circ \st =1\sb F$. Replacing $\Phi $ by $\Phi \sp
{-1}$ and vice versa, we change places $\st $ and $\tau$ and therefore obtain  $\Phi \sp {-1} (\st) \circ \tau
=1\sb F$ and hence $\st \circ \Phi (\tau) =1\sb {\Phi (F)}$. Thus $\st $ is an isomorphism from $F$ onto $\Phi
(F)$ and the statement is true.
\end{proof}

The result above shows that if an automorphisms takes a monogenic free algebra to isomorphic one it does the
same with all free algebras in $\C$.

The next fact is very simple but also very useful.
\begin{lem}\label{subcategory}
Let $\C $ be a subcategory of a category $\D$ and $\Phi :\C \to \D $ be a functor. Let $\E $ be a subcategory of
$\C$ satisfying the following conditions: for every $\E -$algebra $A$, its image $\Phi (A)$ is also an $\E
-$object and there exists a quasi-isomorphism (a $\D -$isomorphism) $\st \sb A :A\to \Phi (A)$ such that $\Phi
(\nu )=\st \sb B \circ \nu \circ \st \sb A \sp {-1}$ for every $\E -$morphism $\nu :A\to B$.

Then $\Phi $ is a composition of two functors $\Phi = \Psi \circ \Gamma $, where $\Gamma :\C \to \D $ is an
identity on $\E $ (preserves all objects and all morphisms of $\E$), and the functor $\Psi$ is an inner
automorphism of the category $\D$.

\end{lem}
\begin{proof}
Under hypotheses, we have a family  $(\st \sb A : A\to \Phi (A)\;\vert A\in Ob\E)$ of $\D -$isomorphisms. We
construct an inner automorphism $\Psi $ of $\D $ in the following way. For every $\D $-object $U$ , we set $\Psi
(U)=U $ if $U$ is not an object of $\E $ , and $\Psi (A) =\Phi (A) $ for every object $A$ of $\E $. Further, we
define $\D$-isomorphism $\tau \sb U :U\to \Psi (U)$ in the following way: $\tau \sb U =1\sb U$ if $U $ is not an
object of  $\E $, and $\tau \sb A = \st \sb A $, where $A\in Ob (\E )$. For every $\D$-morphism $\mu :U\to V$,
let $\Psi (\mu ) = \tau \sb V \circ \mu  \circ \tau \sb U \sp {-1}$. According to the given construction, $\Psi
$ is an inner automorphism of $\D $.

It is clear that $\Phi = \Psi \circ \Psi \sp {-1} \circ \Phi $. Let $\Gamma =\Psi \sp {-1} \circ \Phi $.
According to this definition, we have that $\Gamma (A)=(\Psi \sp {-1} \circ \Phi )(A)=A $ for all $A\in Ob(\E)$
and $\Gamma (\nu )=(\Psi \sp {-1} \circ \Phi )(\nu)=\st \sb B \sp {-1}\circ \Phi (\nu )\circ \st \sb A =\st \sb
B \sp {-1}\circ \st \sb B \circ \nu \circ \st \sb A \sp {-1}\circ \st \sb A =\nu $ for every $\E$-morphism $\nu
:A\to B$.
\end{proof}
The next result is a variation from above.

\begin{lem}\label{subclass} Let  $\Phi $ be an automorphisms of a category $\C$. Suppose further that
for some class $\bf E$ of $\C -$objects,  $\Phi (A)$ is isomorphic to $A$ for every $A\in \bf E$. Then $\Phi$ is
a composition of two $\C -$automorphisms $\Phi = \Psi \circ \Gamma $, where $\Gamma $ leaves fixed all objects
from $\bf E$ and $\Psi$ is an inner automorphism.
\end{lem}
\begin{proof}
If the class $\bf E$ is closed under $\Phi$ and $\Phi \sp {-1}$, apply the previous lemma to the subcategory $\E
$ whose class of objects is $\bf E$ and whose morphisms are identity morphisms only (discrete subcategory) and
obtain the required result where $\D =\C$. If $\bf E$ is not closed we can consider its $\Phi -$ and $\Phi \sp
{-1}-$closure that has clearly the same property.
\end{proof}
\vspace{12pt}

\subsection{Main function}
For all categories $\C$ that are interesting for us, every automorphism $\Phi $ of $\C$ transports a monogenic
free algebra to an isomorphic one. Applying \lemref{subclass}, we restrict our consideration to the case when
the automorphism $\Phi $ under consideration leaves fixed all free algebras in $\C$. This assumption implies
that the maps $s\sb A$ are permutations of all free algebras in $\C$ and the formula~\ref{s-maps}  becomes more
simple:
\begin{equation}\label{maps}
s\sb A (a)= \Phi (\alpha \sb a \sp A )(x\sb 0).
\end{equation}

Of course, the function $A\mapsto s\sb A$ that we call a {\it main function} is not a unique function that
realizes $\Phi$ as a $\D -$inner automorphism for some category $\D$.

\begin{defn}
Let for every $\C-$algebra $A$ we have a permutation $c\sb A$ of its underlying set. The function $A\mapsto c\sb
A $, $ A \in Ob\, \C$ is called {\it central } if for every homomorphism $\nu :A\to B$ the following equation is
satisfied: $c\sb B \circ \nu \circ c \sp {-1} \sb A =\nu$, in other words, if this function determines the
identity automorphism of the category $\C $.
\end{defn}
It is obvious that function $A\mapsto \st \sb A $, $ A \in Ob\, \C$ from Definition~\ref{D-inner} is determined
for a potential inner automorphism $\Phi $ up to a central function. Thus, if we find out that an automorphism
$\Phi$ is $\D -$inner, it may be that it is $\D '-$inner for some subcategory $\D '$ of $\D$ and even inner.
Indeed, being bijections, these maps allow to define a new algebraic structure $A\sp *$ on the underlying set of
every algebra $\Phi (A)$  such that $s\sb A :A \to A\sp *$ is an isomorphism.

\begin{lem}\label{InnerAndCentral}
An automorphism $\Phi$ of $\C$ is $\D -$inner for an extension $\D$ of $\C$ (by adding some kind of
quasi-homomorphisms) if and only if there exists a central function $A\mapsto c\sb A$, $A\in Ob\, \C$ such that
every $c\sb {\Phi (A)}$ is a quasi-isomorphism of $\Phi (A)$ onto $A\sp *$. Particularly, $\Phi$ is inner if and
only if every $c\sb {\Phi (A)}$ in the condition above is a standard isomorphism.
\end{lem}
\begin{proof}
If $\Phi$ is $\D -$inner, $\Phi (\nu )=\st \sb B \circ \nu \circ \st \sb A \sp {-1 }$ for all $\nu :A \to B $,
where  $\st \sb A : A\to \Phi (A)$ are $\D -$isomorphisms for every $A$. Consider $c\sb {\Phi (A)}=s\sb A \circ
\st \sb A \sp {-1}$. Clearly, $c\sb {\Phi (A)}$ is a quasi-isomorphism of $\Phi (A)$ onto $A\sp *$. Since the
considered two functions $A\mapsto s\sb A$ and $A\mapsto \st \sb A$ define the same automorphism, the function
$A \mapsto c\sb A$ is central.

Inversely, if there exists a central function $A\mapsto c\sb A$, $A\in Ob\, \C$, such that every $c\sb {\Phi
(A)}$ is a quasi-isomorphism of $\Phi (A)$ onto $A\sp *$, then we set $\st \sb A=c\sb {\Phi (A)}\sp {-1}  \circ
s\sb A $. Clearly that $\st \sb A$ is a $\D -$isomorphism of $A$ onto $\Phi (A)$, and $\st \sb B \circ \nu \circ
\st \sb A \sp {-1 }= c\sb {\Phi (B)}\sp {-1}  \circ s\sb B  \circ \nu \circ s\sb A \sp {-1}\circ c\sb {\Phi (A)}
=c\sb {\Phi (B)}\sp {-1}  \circ \Phi (\nu )\circ c\sb {\Phi (A)} =\Phi (\nu )$ for all $\nu :A \to B $. Hence
$\Phi $ is $\D -$inner.

\end{proof}
Therefore to describe an automorphism of a given category, the first step is to describe maps $s\sb A$ defined
by \ref{maps} and then to try to find a suitable central function or to prove that it does not exist. The
formula~\ref{maps} can be rewritten in the following useful form:
\begin{equation}\label{E:maps}
\Phi (\alpha )(x\sb 0)=(s\sb A \circ \alpha ) (x\sb 0)
\end{equation}
for every $\alpha :A\sb 0 \to A$.

The following simple fact shows that we can reduce the problem to a case when maps  $s\sb {A}$ have very nice
conditions. If $A$ is a free algebra over a set $X$ we say as usual that $X$ is a basis of $A$.

\begin{thm}\label{saveGenerators} Fix for every free algebra $A$ in $\C$ some its basis $X\sb A$. Let
an automorphism $\Phi $ of $\C$ preserve all free algebras. Then $\Phi$ is a composition of two automorphisms
$\Phi =\Psi \circ \Gamma $, where $\Psi $ is an inner automorphism and $s\sp {\Gamma } \sb A (x)= x $ for all
$x\in X \sb A$ for every free algebra $A$.
\end{thm}
\begin{proof} According to the second part of \thmref{potent-inner}, for every free algebra $A$ we have
an automorphism $\sigma \sb A$ of $A$ such that $\sigma \sb A (x)=s\sb A (x)$ for all $x\in X\sb A$. Now we
ready to apply \lemref{subcategory}, setting $\D =\C$ and $\E$ is equal to the subcategory of $\C$ containing
all free algebras $A$ in $\C$  with set of morphisms $\{\alpha \sb x \sp A ,\; x\in X \sb A \}$ from the $A \sb
0$ to the others ones plus, of course, the identities $1\sb {A\sb 0}$ and $1\sb A$.

Since $\Phi (\a \sb x \sp A )(x\sb 0 )=s\sb A \circ \a \sb x \sp A \circ s\sb {A\sb 0 }\sp {-1}(x\sb 0 )=s\sb A
\circ \a \sb x \sp A (x\sb 0 )=s\sb A (x) =\st \sb A \circ \a \sb x \sp A (x)= \st \sb A \circ \a \sb x \sp A
\circ 1\sb {A\sb 0} (x\sb 0)$, we obtain that $\Phi (\a \sb x \sp A )=\st \sb A \circ \a \sb x \sp A \circ 1\sb
{A\sb 0}$ for all $x\in X\sb A$ and all objects $A$ of the category $\E$. That is, the restriction of $\Phi$ to
$\E$ acts according to conditions of \lemref{subcategory}.

Applying this Lemma, we obtain that $\Phi$ is a composition of two automorphisms  $\Phi =\Psi \circ \Gamma $,
where $\Psi $ is an inner automorphism and $\Gamma (\alpha \sb x \sp A)= \alpha \sb x $ for all $x\in X\sb A$.
The last condition means that $s\sp {\Gamma } \sb A (x)= x $ for all $x\in X \sb A$ for every free algebra $A$.
\end{proof}

\begin{cor}Let $\Phi$ be an automorphism of the category $\C$ leaving fixed a free algebras $A$ over the set
$X$ and $\Phi (\alpha \sb x )= \alpha \sb x $ for all $x\in X$. Let $f:X\to A $. Denote by $\theta \sb f $ the
unique endomorphism of $A$ such that $\theta \sb f (x)=f(x)$ for all $x\in X$. Then $\Phi (\theta \sb f)=\theta
\sb {s\sb A \circ f}$.
\end{cor}
\begin{proof}The condition $\theta \sb f (x)=f(x)$ can be expressed by equality $\alpha \sb {f(x)}= \theta \sb f \circ
\alpha \sb x $. Apply $\Phi$ to this equality and obtain $\alpha \sb {s\sb A (f(x))}= \Phi (\theta \sb f ) \circ
\alpha \sb x$. Hence $s\sb A (f(x))=\Phi (\theta \sb f )(x)$. Since the last one is valid for all $x\in X $, we
have $\Phi (\theta \sb f)=\theta \sb {s\sb A \circ f} $.
\end{proof}
With the preceding notations, let $X=\{ x\sb 1 ,\ldots ,x\sb n \}$, $f(x\sb 1)=a\sb 1 ,\ldots , f(x\sb n )=a\sb
n $. In this situation we write $\theta \sb {a\sb 1 ,\ldots ,a\sb n }$ instead of $\theta \sb f $. Thus we have
the following result:
\begin{equation}
\Phi (\theta \sb {a\sb 1 ,\ldots ,a\sb n })=\theta \sb {s \sb A (a\sb 1 ),\ldots ,s\sb A (a\sb n )}.
\end{equation}
\subsection{Derivative algebras}
From this moment on we consider the category $\Theta \sp 0(\V )$ defined in \secref{Intro}, that is, the full
subcategory of $\Theta (\V )$ formed by all free algebras in $\V$ over finite subsets of a fixed infinite set
$X\sb 0$. This restriction is motivated by the future applications only. According to results obtained in the
previous section, we can restrict our consideration to the case that automorphisms $\Phi $ of this category
satisfy the following two conditions:

1) $\Phi $ leaves fixed all objects ,

2) for every algebra $A$, the its basis $X\sb A \subset X\sb 0$  is fixed and $s\sb A \sp {\Phi}$ leaves fixed
all elements of $X\sb A$.

Let $\Phi$ be an automorphism of such kind. It determines two new structures on underlying set of every algebra
$A$. This structures have the same type that the source structure. First from them was defined above, namely, it
gives the algebra $A\sp * $ induced by permutation $s\sb A $, thus $s\sb A :A\to A\sp *$ is an isomorphism. The
second one we define below.

Let $\omega$ be the symbol of a basic $k-$ary operation. Denote by $\omega \sb A$ the corresponding $k-$ary
operation of the algebra $A$. Consider an algebra $A$ which set $X\sb A =\{x\sb 1 ,\ldots ,x\sb n\}$ of free
generators contains not less elements than $k$. Fix the term $\omega (x\sb 1 ,\ldots ,x\sb k)$ and corresponding
element $w=\omega \sb A (x\sb 1 ,\ldots ,x\sb k)$ in $A$. For every elements $a\sb 1,\ldots ,x\sb k \in A$, we
have:
\begin{equation}\label{value}
\omega \sb A (a\sb 1 ,\ldots ,a\sb k)=\theta \sb {a\sb 1 ,\ldots ,a\sb k ,a\sb {k+1},\ldots ,a\sb n}(\omega \sb
A (x\sb 1 ,\ldots ,x\sb k)),
\end{equation}
where $a\sb {k+1}=\ldots =a\sb n =a\sb k$.

Now apply $\Phi$ and consider the element $\tilde{w}=s\sb A (w)= s\sb A (\omega \sb A (x\sb 1 ,\ldots ,x\sb
k))$. Being an element of the free algebra $A$, it is also a term. And we can define by means of it a new
$k-$ary operation $\tilde{\omega} \sb A$:

\begin{equation}\label{newvalue}
\tilde{\omega }\sb A (a\sb 1 ,\ldots ,a\sb k)=\theta \sb {a\sb 1 ,\ldots ,a\sb k ,a\sb k,\ldots ,a\sb
k}(\tilde{w}).
\end{equation}

Being defined by means of a term, the new operation is a derivative operation which is often called a polynomial
operation. Let $B$ be another algebra. The term $\tilde{w}$ determines a new operation $\tilde{\omega }\sb B $
by the usual way. Let $b\sb 1,\ldots ,b\sb k \in B$. Consider a homomorphism $\nu: A\to B$ defined by $\nu (x\sb
1 )=b\sb 1 ,\ldots ,\nu (x\sb k )=b\sb k ,\nu (x\sb {k+1})=\ldots =\nu (x\sb n )= b\sb k$ and set
\begin{equation}\label{derivativeoper}
\tilde{\omega }\sb B (b\sb 1 ,\ldots ,b\sb k)=\nu (\tilde{\omega }\sb A (x\sb 1 ,\ldots ,x\sb k)).
\end{equation}
It seems at first sight that this operation depends on algebra $A$ we have chosen. But the next result shows
that it is not so.

\begin{thm}\label{main} For every algebra $B$, the derivative operation $\tilde{\omega }\sb B $, defined by
\ref{derivativeoper}, coincides with the induced operation $\omega \sb B \sp *$, that is $s\sb B (\omega \sb B
(b\sb 1 ,\ldots ,b\sb k))=\tilde{\omega }\sb B (s \sb B (b\sb 1 ),\ldots , s\sb B (b\sb k))$ for every $b\sb
1,\ldots ,b\sb k \in B$.
\end{thm}

\begin{proof} The homomorphism $\nu :A\to B $ in \ref{derivativeoper}  is defined by $\nu (x\sb
1 )=b\sb 1 ,\ldots ,\nu (x\sb k )=b\sb k ,\nu (x\sb {k+1})=\ldots =\nu (x\sb n )= b\sb k$. Further,  $\Phi (\nu
)=s\sb B \circ \nu \circ s\sb A \sp {-1}$ and thus $\Phi (\nu ) (x\sb 1) =s\sb B (b\sb 1),\ldots ,\Phi (\nu )
(x\sb k)=s\sb B (b\sb k), \Phi (\nu ) (x\sb {k+1})=\ldots =\Phi (\nu )(x\sb n )=s\sb B ( b\sb k )$. Thus we
obtain
\begin{multline}
s\sb B (\omega \sb B (b\sb 1 ,\ldots ,b\sb k))=s\sb B (\nu (\omega \sb A (x\sb 1 ,\ldots ,x\sb k)))=\\
=\Phi (\nu ) \circ s\sb A (\omega \sb A (x\sb 1 ,\ldots ,x\sb k))=\Phi (\nu )(\tilde{\omega }\sb A (x\sb 1 ,\ldots ,x\sb k))=\\
=\tilde{\omega }\sb B (s \sb B (b\sb 1 ),\ldots , s\sb B (b\sb k)).
\end{multline}
Thus $\omega \sb B \sp * =\tilde{\omega} \sb B $ is a derivative operation on $B$.
\end{proof}

Now we like to explain the fundamental importance of the previous result. For every algebra $A$ of our category
we have a derivative algebra $\tilde{A}$ of the same type. Permutations $s\sb A$ are isomorphisms of $A$ onto
$\tilde{A}$. Therefore we know what quasi-homomorphisms we need to add to make $\Phi$ a $\D -$inner automorphism
for some extension $\D$ of our category. The problem is reduced to finding the derivative operations. To do it
we have to consider a free algebra $A$ which basis $X\sb A$  contains as many elements as the maximum of arities
of all operation is. All derivative operations are defined by terms (polynomials) of the source structure, the
algebra $\tilde{A}$ is isomorphic to $A$ and the isomorphism $s\sb A :A\to \tilde{A}$ preserves all free
generators from $X\sb A$ which of course are free generators for $\tilde{A}$.

Further, every source operation $\omega \sb A$ is a polynomial operation with respect to algebra $\tilde{A}$.
Replace all operations $\tilde{\omega}$ in this polynomial by their expression as polynomials in $A$. We obtain
a word $w$. It means that identities of the kind $\omega (x\sb 1 ,\ldots ,x\sb k )=w$ are satisfied in the
variety under consideration. Studying identities of this kind helps us to find the derivative structure.

Of course, all mentioned above is valid also for every full subcategory of $\Theta (\V )$, containing a
monogenic free algebra and a free algebra which set of free generators contains not less elements than arities
of all basic operations. For the category $\Theta \sp 0(\V )$ specially, we can immediately conclude that every
term $s\sb A (\omega \sb A (x\sb 1 ,\ldots ,x\sb k))$ is build also from the same free generators.

\subsection{Some simple examples}
Now we present three examples only to explain our method because the results we obtain here are known and
published, so we can compare the ways that lead to the same results.

{\it Semigroups}. For the simplest example, consider the variety $\mathbf {SEM}$ of all semigroups. The unique
identity of the mentioned kind satisfied in $\mathbf {SEM}$ with two variables is $xy =xy$. That means we have
only one new derivative operation $x\bullet y =yx$. Thus all automorphisms of the category $\mathbf {SEM}\sp 0$
of free semigroups are $\D$-inner, where the category $\D$ is the category whose objects are as before free
semigroups but whose morphisms are both homomorphisms and anti-homomorphisms of semigroups.

If an automorphism $\Phi$ determines the same derivative operation it is inner. If it determines the dual
operation, it is not inner because in this case there are no central map which is an isomorphism of the
two-generated free semigroup onto the dual semigroup. In the case of variety of commutative semigroups, there
are only inner automorphisms.

{\it Groups}. The similar reasons show that all automorphisms of the category of all free finitely generated
free groups are inner because the map $g\mapsto g\sp {-1}$  is an isomorphism of a group on the dual one, and
the function that assigns to every group such map is a central function.

{\it Lie algebras }. More complicated is the case when $\V$ is the variety of Lie algebras over an infinite
field $\K$ because there are now two binary operations, "+" and "[~]", and the set of unary operations $a\mapsto
ka$ for every $k\in \K$. To obtain the derivative unary operations, we consider monogenic free Lie algebra,
which is one-dimensional linear space $L=\K x$ with trivial multiplication. Thus $s\sb L (kx) =\vp (k)x$. The
map $\vp :\K \to \K$ is of course a bijection preserving multiplication in $\K$. Hence $\vp (0)=0$ and $\vp
(1)=1$.

To find derivative binary operations, we consider two-generated free Lie algebra $F=F(x,y)$. We have $x\bot
y=s\sb F (x+y)=P(x,y)$, where $P(x,y)$ is a polynomial in $F$. Since $x\bot y $ is homogeneous of the degree 1,
the polynomial $P(x,y)$ is linear: $x\bot y =ax+by$. It is clear that $a=b=1$ because of commutativity and of
the condition $x\bot 0=x$.

Conclusion: $s\sb F $ is an additive isomorphism, such that $s\sb F (kw)=\vp (k)s\sb F (w)$, where $\vp$ is an
automorphism of $\K$, $k\in \K$ and $w\in F$.

Further, since $x\ast y =s\sb F ([xy])$ is a homogeneous polynomial of degree 2, $x\ast y =a[xy]$, where $a\in
\K$. Consider a function which assigns to every free Lie algebra $F(X)$ the permutation ${\bf c}\sb {F(X)}$ of
$F(X)$ as follows: ${\bf c}\sb {F(X)}(w)=w/a$. This function is clearly central. Hence every automorphism of the
category $\Th \sp 0(\V)$ is semi-inner according to definition given in \cite{LieAlg}, that is in
Definition~\ref{D-inner}, the bijections $\st \sb A$ are additive and multiplicative isomorphisms of the
corresponding Lie rings but $\st \sb A (kw)=\vp (k)\st \sb A (w)$ for every $k\in \K$ and $w\in A$.

In the next two sections, we apply the method to the variety of all associative linear algebras over a fixed
infinite field and to the variety of all representations of groups in unital $R-$modules over a associative
commutative ring $R$ with unit.

\section{Associative linear algebras}\label{secasslinalg}
In this section,  $\K$ will always denote an infinite field. Let $A$ be an associative ring with unit and $f:\K
\to A$ be a ring homomorphism of $\K $ into the center of the ring $A$. We consider such a homomorphism as an
associative (linear) algebra over $\K$, shortly, a $\K -$algebra. Given two $\K -$algebras $f:\K \to A$ and
$g:\K \to B$. A homomorphism of the first algebra into second one is a ring homomorphism  $\nu :A\to B $ such
that $g=\nu \circ f$. As usual, we identify elements of the field $\K$ considered as symbols of 0-ary operations
(constants) with their images in a ring and write $ka$ instead of $f(k)a$, if it does not lead to a
misunderstanding.

The category of all $\K -$algebras and their homomorphisms we denote by $Ass-\K$ . Let  $(Ass-\K )\sp 0$ be the
full subcategory of $Ass-\K$ consisting of all free $\K -$algebras over finite subsets of an infinite set $X\sb
0$. We fix in $(Ass-\K )\sp 0$ a monogenic free algebra $F\sb 1$ with a free generator $x\sb 0$ and a
two-generated free algebra $F\sb 2$ whose two free generators are $x$ and $y$.

Let $\Phi$ be an automorphism of the category $(Ass-\K )\sp 0$. It is clear that $\Phi (F\sb 1 )$ is isomorphic
to $F\sb 1 $, hence $\Phi $ is potential inner (\thmref{potent-inner}) and we can assume that $\Phi $ leaves
fixed all objects of the category and the permutations $s\sb A$ act as identities on bases $X\sb A$ of objects
$A$ of this category. Every map $s\sb A$ is an isomorphism of $A$ onto the derivative algebra $\tilde{A}$.
According to the method suggested in the previous section, we have to find the derivative algebraic structure
$F\sb 2 \sp *=\tilde{F\sb 2}$, which is of course an associative free two-generated $\K-$algebra isomorphic to
$F\sb 2$.
\begin{lem}\label{ActsOnField} Every map $s\sb A$ acts on  $\bar{\K} =f\sb A (\K )$ as a permutation.
\end{lem}
\begin{proof} Let $u\in A$. It is obvious that $u\in \bar{\K} $ if and only  if for every
endomorphism $\nu $ of $A$ it holds: $\nu \circ \alpha \sb u =\alpha \sb u$. Apply automorphism $\Phi $ and
obtain that for every endomorphism $\nu $ of $A$ it holds: $\Phi (\nu )\circ \alpha \sb {s\sb A (u)} =\alpha \sb
{s\sb A (u)}$. Since $\Phi (\nu )$ runs all endomorphisms, it means that $u\in \bar{\K} $ if and only if $s\sb
A(u) \in \bar{\K} $.
\end{proof}
This result shows that the automorphism $\Phi $ determines a permutation $ \tilde{s}$ of the field $\K$. The
derivative 0-ary structure on the set $A$ is given by the map $f\sb A \sp * : \K \to A \sp *$, where $f\sb A \sp
* =s\sb A \circ f\sb A =f\sb A\circ \tilde{s}$. We denote $\tilde{0}=\tilde{s} (0)$ and $\tilde{1}=\tilde{s}
(1)$. We remind that the same symbols will denote the corresponding elements in $\K -$algebras.

\vspace{12pt} {\bf Remark}. The same result can be obtained using the following reasons. The free algebra in
$Ass-\K$ over the empty set of generators is $\K $. Thus $s\sb {\K} $ is a permutation of $\K $ and according to
the general definition \ref{derivativeoper} of derived operations, $f\sb A \sp * =\tilde{f} \sb A=f\sb A \circ
s\sb {\K} = s\sb A \circ f\sb A$. Here $s\sb {\K} $ coincides with $ \tilde{s}$ above.

\vspace{12pt} Denote for simple $s\sb {F\sb 2}$ by $s$. Now our aim is to find derivative operations:  $x\bot
y=s(x+y)$ and $x\odot y=s(xy)$, that are polynomials of two non-commutative variables $x,y$. Since $F\sb 2 \sp
*$ is a free $\K-$algebra free generated by $x$ and $y$, the terms $x+y$ and $xy$ are also terms (polynomials)
in $F\sb 2 \sp *$. Repeat for this case the reasonings in the previous section. Namely, $x+y$ is a polynomial
$w\sp *$ constructed by means of $\oplus $ and $\odot $. Replacing this operations by their expressions as
polynomials in $F\sb 2 $, we obtain a polynomial $w$ degree of which is not less than degree of $w\sp *$. Since
we have an identity $x+y =w$, the degree of $w$ is equal to 1. Therefore the degree of $x\bot y$ is not greater
than 1. The same reasons lead to the conclusion that the degree of $x\odot y $ is not greater than 2. It is a
crucial moment in our considerations.

Consider $x\bot y=ax+by+c$, where $a,b,c \in \bar{\K}$.  We obtain that $ax+b\tilde{0}+c=x$ and $a\tilde{0} +by
+c =y$. These identities give that $a=b=1$ and $c=-\tilde{0}$. Thus $x\bot y= x+y-\tilde{0} $.

Now consider $x\odot y=s(xy)= a\sb {11}x\sp2+a\sb {12}xy +a\sb {21}yx +a\sb {22}y\sp 2 +a\sb 1 x +a\sb 2 y +a$.
Since $x\odot \tilde{0}=\tilde{0}\odot y=\tilde{0}$, we have
\begin{multline}\label{zero}
a\sb {11}x\sp 2 +a\sb {12} x\tilde{0} + a\sb {21} \tilde{0} x +a\sb {22}(\tilde{0})\sp 2 + a\sb 1 x +a\sb 2
\tilde{0} +a =\\
=a\sb {11}(\tilde{0})\sp 2 +a\sb {12} \tilde{0} y + a\sb {21} y\tilde{0} +a\sb {22}y\sp 2 + a\sb 1 \tilde{0}
+a\sb 2 y +a =\tilde{0}
\end{multline}
\noindent and hence $a\sb {11}=a\sb {22}=0,\; a\sb {12}\tilde{0} +a\sb {21} \tilde{0} +a\sb 1 =0, \; a\sb
{12}\tilde{0} +a\sb {21} \tilde{0} +a\sb 2 =0,\; a\sb 1 \tilde{0} +a=\tilde{0}$. And as a consequence from the
last equations, we have $a\sb 1 =a\sb 2=b$.

It gives the equality: $x\odot y= a\sb {12}xy +a\sb {21}yx +b (x+y) +a$.

Now we use the associativity law. Observe that $x\odot x=a\sb {12}x\sp 2 +a\sb {21}x\sp 2+ 2b x +a$. The
identity $(x\odot x)\odot y =x\odot (x\odot y)$  gives:
\begin{multline}\label{ass}
a\sb {12}(a\sb {12}x\sp 2 +a\sb {21}x\sp 2+ 2b x +a)y+a\sb {21}y (a\sb {12}x\sp 2 +a\sb {21}x\sp 2+ 2b x +a)+\\
+b (a\sb {12}x\sp 2 + a\sb {21}x\sp 2+ 2b x +a +y)+a= \\=a\sb {12}x (a\sb {12}xy +a\sb {21}yx +b (x+y) +a)+ a\sb
{21}(a\sb {12}xy +a\sb {21}yx +b (x+y) +a)x + \\+b(x+a\sb {12}xy +a\sb {21}yx +b (x+y) +a)+a.
\end{multline}
Comparing coefficients by  $x\sp 2y$, we obtain $a\sb {12} \sp 2 + a\sb {12}a\sb {21}=a\sb {12} \sp 2$ and hence
$a\sb {12}a\sb {21}=0$. Thus one and only one of the coefficients $a\sb {12}$ and $a\sb {21}$ vanishes. Suppose
it is $a\sb {21}$. In this case, $x\odot y=kxy +bx +by+a$ and $k\tilde{0} +b =0, \;b\tilde{0}+a=\tilde{0}$. We
have $x\odot y=kxy -k\tilde{0} x -k\tilde{0} y+\tilde{0}+k(\tilde{0})\sp 2 +\tilde{0}=
k(x-\tilde{0})(y-\tilde{0})+\tilde{0}$. It is clear that $k=(s(1)-s(0))\sp {-1} \not =0$. In the case $a\sb {12}
=0$, the similar conclusion takes place.

Finally, we obtain
\begin{equation}\label{derivring}
x\bot y= x+y -\tilde{0}\quad and \quad x\odot y= k(x-\tilde{0})(y-\tilde{0})+\tilde{0}
\end{equation}
or in dual case for multiplication
\begin{equation}\label{dualderivring}
x\odot y= k(y-\tilde{0})(x-\tilde{0})+\tilde{0}.
\end{equation}
{\bf Conclusion:} We have found the derivative structure $F\sb 2 \sp*$ defined by means of the permutation $s$,
its operations are defined above. It is interesting to mention that the derived operations are obtained from the
source operations by means of translation and contraction operators. As a result, we have found the main
function. This function assigns to every algebra $A$ a permutation $s\sb A$  satisfying the following
conditions:
$$s\sb A \circ f\sb A  =f\sb A \circ \tilde{s},$$
$$s\sb A (u+v)=s\sb A (u)+s\sb A (v)-\tilde{0},$$
$$s\sb A (uv)=(\tilde{1}-\tilde{0})\sp {-1}(s\sb A (u)-\tilde{0})(s\sb A (v)-\tilde{0})+\tilde{0}$$
or for dual case $$s\sb A (uv)=(\tilde{1}-\tilde{0})\sp {-1}(s\sb A (v)-\tilde{0})(s\sb A
(u)-\tilde{0})+\tilde{0}.$$

Adding such kind of maps to the category $(Ass-\K )\sp 0 $, we obtain an extension $\Q $ of $(Ass-\K )\sp 0 $
and can formulate the first description:
\begin{lem} All automorphism of the category $(Ass-\K )\sp 0$ are $\Q -$inner.
\end{lem}

{\bf Remark }. "Such kind of maps" means that new maps $\st : A\to B $ (quasi-homomorphisms) satisfy the
following conditions:
$$\st ( f\sb A  (\K ))=f\sb B (\K ),$$
$$\st(u+v)=\st(u)+\st(v)-\st (0),$$
$$\st (uv)=(\st (1)-\st (0))\sp {-1}(\st (u)-\st (0))(\st (v)-\st (0))+\st (0)$$
or $$\st (uv)=(\st (1)-\st (0))\sp {-1}(\st (v)-\st (0))(\st (u)-\st (0))+\st (0).$$

It is easy to see, that in the case $\tilde{0}=0,\tilde{1}=1$, all maps $s\sb A$ are isomorphisms or all are
anti-isomorphisms. In this case the automorphism $\Phi$ is inner or induced by means of anti-isomorphisms.
Consider the second case. It is known that there exists a special automorphism of our category named the mirror
automorphism. We repeat its construction.

Given a free $\K -$algebra $A$ with the set $X\sb A$ of free generators. Consider free semigroup $X\sb A \sp +$.
Correspond to every word $w=x\sb {i\sb 1}\ldots x\sb {i\sb n}$ the {\it inverse } word $\bar{w}=x\sb {i\sb
n}\ldots x\sb {i \sb 1}$, that is, all letters are written in inverse order. The map $w\mapsto \bar{w}$ is an
anti-automorphism of this free semigroup. This map can be uniquely extended to an anti-automorphism $\eta \sb A
$ of the $\K -$algebra $A$. It is obvious that for every homomorphism $\nu :A\to B$ of free algebras, the map
$\eta \sb B \circ \nu \circ \eta \sb A \sp {-1}$ is also an homomorphism from $A$ into $B$. Hence we have an
automorphism $\Upsilon $ of the category $(Ass-\K )\sp 0$, defined by: $\Upsilon (A)=A$ and $\Upsilon (\nu
)=\eta \sb B \circ \nu \circ \eta \sb A \sp {-1}$ for every $\nu :A\to B$. This automorphism $\Upsilon $ is
called the {\it mirror} automorphism. Take note of the fact that all $\eta \sb A$ leave fixed all elements from
$X\sb A$.

The automorphism $\Gamma =\Phi \circ \Upsilon  $ satisfies the same conditions that we assume for $\Phi $ but
$s\sb A \sp {\Gamma}$ are homomorphisms onto derived algebras with operations~\ref{derivring}. Since $\Phi
=\Gamma \circ \Upsilon $, we can concentrate our attention on the case when the maps $s\sb A $ are isomorphisms
onto such derived algebras.

Assign to every $\K -$algebra $A$ the permutation $c \sb A$ of $A$ defined as follows:
$$c\sb A (u)=\frac {u}{\tilde{1}-\tilde{0}} +\tilde{0}
$$
for every $u\in A$.
\begin{lem}\label{centralmap}
Every map $c\sb A$ is an ring-isomorphism of  $A$ onto $\tilde{A}$, and the function $A\mapsto c\sb A$  is a
central function.
\end{lem}
\begin{proof}
Indeed, denote for shot $k=\tilde{1}-\tilde{0}$ and calculate:
\begin{multline}
c\sb A (u+v)=\frac {u+v}{k} +\tilde{0}=\frac {u}{k}+\tilde{0} + \frac {v}{k} +\tilde{0}-\tilde{0}=\\
= c\sb A (u)+c\sb A (v)-\tilde{0}=c\sb A (u)\bot c\sb A (v).
\end{multline}
Further,
\begin{multline}
c\sb A (uv)=\frac {uv}{k} +\tilde{0}=k(\frac {u}{k}+\tilde{0} -\tilde{0})(\frac {v}{k}+\tilde{0}
-\tilde{0})+\tilde{0}=\\=k(c\sb A (u)-\tilde{0})(c\sb A (v)-\tilde{0})+\tilde{0}=c\sb A (u)\odot c\sb A (v).
\end{multline}
Let $\nu :A\to B$ be a homomorphism of $\K -$algebras. For every $w\in A $, we have: $c\sb B \circ \nu \circ c
\sb A \sp {-1} (w)=c\sb B \circ \nu (kw-k\tilde{0})=c\sb B (k\nu (w)-k\tilde{0})=\nu (w)$. Thus $c\sb B \circ
\nu \circ c\sb A \sp {-1}=\nu$ and the considered function is central.
\end{proof}
However, the maps $c\sb A$ in general are not isomorphisms of $\K -$algebras, because
$$c\sb A \circ f\sb A =f\sb A\circ c\sb {\K}=\tilde{f}\sb A  \circ \tilde{s} \sp {-1}\circ c\sb {\K}.$$
Since $c\sb {\K} :\K \to \tilde {\K}$ and $\tilde{s}:\K \to \tilde {\K}$ are isomorphisms, the map $\varphi
=\tilde{s} \sp {-1}\circ c\sb {\K}$ is an automorphism of the field $\K$. Thus $c\sb A \circ f\sb A
=\tilde{f}\sb A \circ \varphi $. This conditions lead to the following definition.
\begin{defn}Let $f\sb A :\K \to A$ and $f\sb B :\K \to B$ be $\K -$algebras. A map $\st : A\to B$ is
said to be a twisted homomorphism if it is a ring homomorphism of $A$ into $B$ and $\st  \circ f\sb A =f\sb B
\circ \varphi $ for some automorphism $\varphi$ of the field $\K$, in other words, $\st (aw)=\varphi (a) \st
(w)$ for every $a\in \K$ and $w\in A$. In this case we say that $\st $ is a $\varphi -$homomorphism.
\end{defn}
Build now the category, denoted by $R(Ass-\K )\sp 0$, that is an extension of the source category obtained by
adding all twisted homomorphisms. According to \lemref{InnerAndCentral}, our automorphism $\Phi$ is an $R(Ass-\K
)\sp 0 -$inner automorphism. More exactly, for every homomorphism $\nu :A\to B$ of $\K -$algebras, the action of
$\Phi $ can be expressed as follows:
\begin{equation}\label{finaldescr}
\Phi (\nu )=\tau \sb B \circ \nu \circ \tau \sb A \sp {-1},
\end{equation}
here $\tau \sb C =s\sb C \sp {-1} \circ c\sb C$ is a $\varphi -$automorphism of the $\K -$algebra $C$ and
automorphism $\varphi $ of the field $\K $ defined by  $ \varphi (a)=\tilde{s}\sp {-1}(\frac{a}{k}+\tilde{0})$
for $a\in \K $. An automorphism $\Phi$ of the category $R(Ass-\K )\sp 0$ defined above is called semi-inner.
Thus we have final description.
\begin{thm}\label{FreeAssAlg} Every automorphism $\Phi$ of $(Ass-\K)\sp 0$ is semi-inner (in particular, inner) or
a composition of a semi-inner and the mirror automorphism.
\end{thm}
It is useful to give the description above in more convenient form. Let $\Phi $ is semi-inner automorphism in
the form \ref{finaldescr} where all $\tau \sb C $ are $\varphi -$automorphisms. The automorphism $\varphi $ of
the field $\K$ can be uniquely extended to a twisted isomorphism $\varphi \sb C$ for every free $\K -$algebra
$C$ by identity acting on $X\sb C \sp +$. These twisted isomorphisms determines a semi-inner automorphism
$\hat{\varphi}$ called standard $\varphi -$automorphism of the category $(Ass-\K)\sp 0$. It is obvious that the
composition $\Phi \circ \hat{\varphi}\sp {-1}$ is an inner automorphism of $(Ass-\K)\sp 0$. Thus we can assert
that
\begin{thm}\label{AutFreeAssAlg} Every automorphism $\Phi $ of the category $(Ass-\K)\sp 0$ can be represented as a
composition of three automorphisms:
$$ \Phi =\Upsilon \circ \hat{\varphi}\circ \Psi  ,$$
where $\Psi $ is inner, $\hat{\varphi}$ is the standard $\varphi -$automorphism and $\Upsilon $ is the mirror or
the identity automorphism.
\end{thm}
\section{Category of group representations}\label{secGroupRep}

\subsection{Basic definitions}

$R$ will always denote associative and commutative ring with unit 1. All $R-$modules under consideration are
supposed to be unital. A representation of a group $G$ in an $R-$module $A$  is an arbitrary group homomorphism
$\rho: G\to Aut\sb R (A)$, where $Aut\sb R (A)$ is the group of all $R-$module automorphisms of $A$.

If such a representation is given we have an {\it action} of the group $G$ in $A$, that is, a map $(a,g)\mapsto
a\cdot g $ from $A\times G $ into $A$, satisfying the following conditions:
\begin{enumerate}
\item for every $g\in G$, the map $a\mapsto a\cdot g$ is an automorphism of the module $A$;

\item $(a\cdot g\sb 1 )\cdot g\sb 2 =a\cdot (g\sb 1 g\sb 2 )$ for every $g\sb 1 , g\sb 2 \in G$ and $a\in A$.
\end{enumerate}

We consider a representation as a triple $(A, G,\cdot )$, where $A$ is an $R-$module, $G$ is a group and $\cdot
$ is an action of $G$ in $A$. Thus a representation can be regarded as an universal algebra $A$ with one binary
operation $+$ and two families of unary operations: actions of elements of the ring $R$ and actions of elements
of a group $G$. However we will regard homomorphisms of representations not as homomorphisms of such universal
algebras but (according to \cite{VarRep}) as homomorphisms of two-sorted algebras.

Given two representations $(A,G,\cdot )$ and $(B,H,\bullet )$. A homomorphism $\mu :(A,G,\cdot )\to (B,H,\bullet
)$ is a pair of maps $(\mu \sp {(1)},\mu \sp {(2)})$ of the kind $\mu \sp {(1)}:A\to B ,\;\mu \sp {(2)}:G \to H
$, satisfying the following conditions:
\begin{enumerate}
\item $\mu \sp {(1)}: A\to B $ is an homomorphism of the $R-$modules ;

\item $\mu \sp {(2)}: G\to H$ is an homomorphism of the groups;

\item $\mu \sp {(1)}$ and $\mu \sp {(2)}$ are connected in the following way: for every $a\in A $ and every
$g\in G $, $$ \mu \sp {(1)}(a\cdot g) =\mu \sp {(1)}(a)\bullet \mu \sp {(2)}(g). $$ In other words the following
diagram is commutative:
$$
\CD
A @>\tilde{g}>>    A\\
@V\mu \sp {(1)} VV     @VV \mu \sp {(1)}V\\
B @>\widetilde{\mu \sp {(2)}(g)}>>    B,
\endCD
$$
where $\tilde{g}$ and $\widetilde{\mu \sp {(2)}(g)} $ denote the corresponding actions.

\end{enumerate}

The class of all group representations as objects and the defined above homomorphisms as morphisms form a
category which we denote by $Rep-R$. The forgetful functor is usually  defined as for two-sorted theory. It
means that it assigns to every representation $(A,G,\cdot )$ the pair of sets $(A,G)$ and to every homomorphism
$(\mu \sp {(1)},\mu \sp {(2)})$ the pair of corresponding maps. Therefore we have a notion of free objects
$(W,F, \star)$ in this category, so called free representations over any pair of sets $(Y,X)$. It means that the
set $Y$ generates $RF-$module $W$, where $RF$ is the group algebra over the ring $R$, the set $X$ generates the
group $G$ and for every representation $(A,G,\cdot )$ and every two maps $f\sp {(1)}: Y\to W$ and $f\sp
{(2)}:X\to F$, there exists a unique homomorphism $\mu :(W,F, \star) )\to (A,G,\cdot )$ such that $\mu \sp
{(1)}$ extends $f\sp {(1)}$ and $\mu \sp {(2)}$ extends $f\sp {(2)}$. But the mentioned functor is not a functor
to the category of sets and maps. Thus we have a problem that we can not apply results from \secref{basic}
immediately.

We introduce the following forgetful functor: it assigns to every object $(A, G,\cdot )$ the set $A\times G$ and
to every homomorphism $\mu :(A,G,\cdot )\to (B,H,\bullet )$ the map $\vert \mu \vert =\mu \sp {(1)}\times \mu
\sp {(2)}:A\times G \to B\times H $. We see that free objects with respect to such forgetful functor do not
coincide with free objects with respect to the two-sorted theory. It leads to some modifications in definitions.
Therefore we check in some cases if the results obtained in \secref{basic} are valid.

In the same way as in the previous section, we consider the category of $(Rep-R )\sp 0$ of all free
representations over pairs $(Y,X)$ of finite subsets of fixed infinite sets $Y\sb 0 $ and $X\sb 0$ respectively.
Let $(W\sb 1, F\sb 1)$ denote a monogenic free representation. That means, that $F\sb 1$ is the infinite cyclic
group $ \{x\sp n \vert n\in \Z\}$,  $W\sb 1 =RF\sb 1$ can be identified with the group algebra over the ring $R$
and the action of this group on the group algebra is the group algebra multiplication. Denote $e=x\sp 0$ the
unit of the group $F\sb 1$ and $1\sb R $ the unit if the ring $R$. Identify for every $r\in R$,  $re$ with $r$,
thus $R$ is embedded in $W\sb 1 $. The pair of singular sets  $(\{1\sb R \},\{ x\})$ is a basis of the
representation $(W\sb 1, F\sb 1)$.

Given a group representation $(A,G,\cdot)$  which we for short denote by $AG$. According to our usual
definitions, denote by $\a \sp {AG} \sb {(a,g)}$ the unique homomorphism from $(W\sb 1, F\sb 1)$ to
$(A,G,\cdot)$ that takes $1\sb R$ to $a\in A$ and $x$ to $g\in G$.

It was mentioned that we consider the product $A\times G$ as the underlying set of a presentation $(A,G,\cdot
)$.  Since forgetful functor is presented by $(W\sb 1, F\sb 1)$, the first part of \thmref{potent-inner} is
valid. It means an automorphism $\Phi $ of the category $(Rep-R )\sp 0$ is potential inner if and only if it
takes the representation $(W\sb 1, F\sb 1)$ to an isomorphic one. The second part of this theorem will be
considered later.

Let $\Phi $ be an automorphism of the category $(Rep-R )\sp 0$ leaving fixed $(W\sb 1, F\sb 1)$. Thus the main
function $(A,G,\cdot) \mapsto s\sb {AG}$ is defined in the same way as in \secref{basic}:
\begin{equation}
s\sb {AG} (a,g)=(\bar{a},\bar{g})\Leftrightarrow \Phi (\a \sp {AG} \sb {(a,g)})=\a \sp {\Phi (AG)} \sb
{(\bar{a}\bar{g})},
\end{equation}
and we have the usual form of $\Phi$ action, namely, for every homomorphisms $\mu :(A,G,\cdot)\to
(B,H,\bullet)$:
\begin{equation}\label{actionRep}
\Phi (\mu )=s\sb {BH} \circ \mu \circ s\sb {AG}\sp {-1}.
\end{equation}

Consider a free representation $(W,F,\cdot )$ with basis $(Y,X)$. Let $f\sp {(1)}:Y\to W $ and $f\sp {(2)}:X\to
G$. The unique endomorphism extending these maps we denote by $\theta \sb {(f\sp {(1)},f\sp {(2)})}$. Hence
$$ (\forall y\in Y ,x\in X )\;\theta \sb {(f\sp {(1)},f\sp {(2)})}(y,x)=(f\sp {(1)}(y),f\sp {(2)}(x)).$$
To apply the method suggested in \secref{basic} we have first to investigate the monoid $End\sb1 $ of
endomorphisms of the free monogenic representation $(W\sb 1 , F\sb 1)$.

\subsection{Some invariants of category automorphisms}
So our first step is studying the monoid $End\sb1 $ of endomorphisms of the free monogenic representation $(W\sb
1 , F\sb 1)$. We denote elements of this monoid by $\nu \sb {(w,g)}$, where $w=\nu \sb {(w,g)}\sp {(1)}(1\sb R)$
and $g=\nu \sb {(w,g)}\sp {(2)}(x)$.

\begin{lem}\label{ideals} The endomorphism $\nu \sb {(0,e)}$, where $0$ is the zero of the ring $R$, is the zero
element of the monoid $End\sb1 $. The endomorphism $\nu \sb {(1\sb R,x)}$ is the unit of this monoid. The set
$T\sb e$ of all endomorphisms $\nu \sb {(w,e)}$, where $w\in W\sb 1$ is a minimal prime ideal in this monoid.
Every prime ideal different from $T\sb e$ contains $T\sb e $ or the set $T\sb 0 =\{\nu \sb {(0,g)}\vert \, g\in
F\sb 1\}$ that also is an ideal of $End\sb1 $.
\end{lem}
\begin{proof} The first and second statements are obvious. Also it is obvious that the sets $T\sb e$ and $T\sb 0$
are ideals in $End\sb1 $. Suppose that for some $\nu ,\mu \in End\sb1 $ we have $\nu \circ \mu \in T\sb e $. It
means that $\nu \sp {(2)}\circ \mu \sp {(2)}(x)=e$ and therefore $\mu \sp {(2)}(x)=e$ or $\nu \sp {(2)}(x)=e$.
Hence $T\sb e $ is a prime ideal. Suppose now that $I$ is another prime ideal and there exists $w\in W\sb 1$
such that $\nu \sb {(w,e)}\notin I$. But $\nu \sb {(w,e)}\circ \nu \sb {(0,g)}=\nu \sb {(0,e)}$ for all $g\in
F\sb 1$. Hence $\nu \sb {(0,g)}\in I$ for all $g\in F\sb 1$, that is, $T\sb 0 \subseteq I$. The conclusion is
that $T\sb e \subseteq I$ or $T\sb 0 \subseteq I$. Since $T\sb 0 \cap T\sb e =\{\nu \sb {(0,e)}\}$, $T\sb e $ is
a minimal prime ideal.
\end{proof}

\begin{lem}\label{unit} $T\sb e$ is a unique minimal prime ideal in $End\sb1 $ that contains many right units.

\end{lem}
\begin{proof} Consider endomorphisms $\nu \sb {(w,e)}$, where $w$ is an element of group algebra $W\sb 1$
satisfying the condition that its value by $x=e$ is e. Since $w=\sum r\sb n x\sp n$, where $n\in \Z $ and $r\sb
n =0$ for almost all $n$, this conditions means that $\sum r\sb n =1$. For every $\nu \in T\sb e$ we have:
$$\nu \sp {(1)} \circ \nu \sb {(w,e)}\sp {(1)}(1\sb R )=\nu \sp {(1)} (w)=\nu \sp {(1)} (1\sb R w)=\nu \sp {(1)} (1\sb
R)\nu \sp {(2)} (w)=\nu \sp {(1)} (1\sb R).$$

That gives us $\nu  \circ \nu \sb {(w,e)}=\nu$, i.e., $\nu \sb {(w,e)}$ is a right unit in $T\sb e$. Suppose
that $I$ is another minimal prime ideal and $\nu \sb {(u,g)}$ is its right unit. According to \lemref{ideals},
$T\sb 0\subseteq I$. Hence $\nu \sb {(0,x)}\in I$ and we obtain: $\nu \sb {(0,x)}=\nu \sb {(0,x)}  \circ \nu \sb
{(u,g)}=\nu \sb {(0,g)}$. This equation gives $g=x$. Thus, our right unit has a form: $\nu \sb {(u,x)}$. Since
$\nu \sb {(u,x)}\circ \nu \sb {(v,x)}=\nu \sb {(uv,x)}$ and  $uv=vu $ for every two elements  $u,v \in W\sb 1$,
every two right units in $I$ coincide.
\end{proof}
Denote by $T\sb x$ the set of all endomorphisms of $(W\sb 1 ,F\sb 1)$ of the kind $\nu \sb {(w,x)}$, where $w$
an arbitrary element of $W\sb 1 $.

{\bf Remark.} The ideal $T\sb 0$ can be determined in the considered monoid because $T\sb 0 =\{\nu \,\vert\,
(\forall \mu \in T\sb e )\,\nu \circ \mu =\mu \circ \nu =\nu \sb {(0,e)}\}$. The set $T\sb x$ can be described
by the following way: $\nu \in T\sb x \Leftrightarrow (\forall \mu \in T\sb 0) \nu \circ \mu =\mu $. Element
$\nu \sb {(0,x)}$ is an unique common element of $T\sb 0 $ and $T\sb x$.

\begin{cor}\label{Tx}
$T\sb x $ is a submonoid of $END\sb 1$ which is multiplicatively isomorphic to the group algebra $RF\sb 1$.
\end{cor}
\begin{proof} It is clear that $T\sb x $ is subsemigroup of $END\sb 1$, containing the unit $\nu \sb {(1\sb R
,x)}$. For every $u,v\in W\sb 1 $ we have:
$$\nu \sb {(v,x)}\sp {(1)} \circ \nu \sb {(v,x)}\sp {(1)}(1\sb R)=\nu \sb {(v,x)}\sp {(1)}(u)=vu,$$
therefore $\nu \sb {(v,x)} \circ \nu \sb {(v,x)}=\nu \sb {(vu,x)}$.
\end{proof}

The next results suppose that we deal with such automorphisms which preserve the object $(W\sb 1 , F\sb 1)$.
\begin{cor}\label{x,eand0} Let an automorphism $\Phi$ of the category $(Rep-R )\sp 0$ leave fixed the object
$(W\sb 1 , F\sb 1)$. Then $\Phi$ preserves the sets $T\sb e$,  $T\sb 0$ and $T\sb x$ and hence the endomorphism
$\nu \sb {(0,x)}$.
\end{cor}
\begin{proof} $\Phi$ acts on the monoid $END\sb 1$ as an automorphism. It preserves $T\sb e$ because of \lemref{unit}.
It preserves the sets $T\sb 0$ and $T\sb x$ and hence the endomorphism $\nu \sb {(0,x)}$ because of Remark
above.
\end{proof}
\begin{cor}\label{MultRingMon} Under hypotheses of the previous result, $\Phi$ induces an multiplicative automorphism
of the group algebra $RF\sb 1$ and hence an automorphism $\varphi $ of the multiplicative monoid of the ring
$R$.
\end{cor}
\begin{proof}  Using \corref{Tx}, define $\varphi (r)=t \Leftrightarrow \Phi
(\nu \sb {(r,x)}) =\nu \sb {(t,x)}$. Then $\varphi (rt)=\varphi (r)\varphi (t)$.
\end{proof}

\subsection{The main function}
In this section we suppose that an automorphism $\Phi$ leaves fixed the monogenic free representation $(W\sb 1 ,
F\sb 1)$. Thus we have the formula~\ref{actionRep} and we start to study the maps $s\sb {AG}$.

Let $AG=(A,G,\cdot)$, $0$ be the zero of the module $A$ and $e$ be the unite of the group $G$. Denote by $T\sp
{AG}\sb e$ and by  $T\sp {AG}\sb 0$ the sets of all homomorphisms of $(W\sb 1 , F\sb 1)$ into $(A,G,\cdot)$ of
the kind $\a \sb {(w,e)}$ and $\a \sb {(0,g)}$ respectively, where $w\in A,\;g\in G$.
\begin{lem} Let $\Phi$ take $AG=(A,G,\cdot)$ to $BH=(B,H,\odot)$. Then $\Phi (T\sp {AG}\sb e)=T\sp {BH}\sb e $
and $\Phi (T\sp {AG}\sb 0)=T\sp {BH}\sb 0 $.
\end{lem}
\begin{proof} Let $\frak M$ be the set of all homomorphisms from $(W\sb 1 , F\sb 1)$ into $AG=(A,G,\cdot)$. Then
$\frak N= \Phi (\frak M)$ is the set of all homomorphisms from $(W\sb 1 , F\sb 1)$ into $BH=(B,H,\odot)$. Since
$T\sp {AG}\sb e =\frak M \circ T\sb e $ and $T\sp {AG}\sb 0 =\frak M \circ T\sb 0 $, we obtain using
\corref{x,eand0} that $\Phi (T\sp {AG}\sb e )=\frak N \circ T\sb e =T\sp {BH}\sb e $ and $\Phi (T\sp {AG}\sb 0)
=\frak N \circ T\sb 0 =T\sp {BH}\sb 0  $.
\end{proof}
The map $s\sb {AG}$ is a subdirect product of two maps $s\sb {AG}\sp {(1)}$ and $s\sb {AG}\sp {(2)}$, first of
which gives the first element of the $s\sb {AG}-$image, and the second map gives the second one: $s\sb
{AG}(u,g)=(\bar{u},\bar{g})\Leftrightarrow s\sb {AG}\sp {(1)}(u,g)=\bar{u}\; \& \; s\sb {AG}\sp
{(2)}(u,g)=\bar{g}$. The both maps are two-place functions. It turns out that we can replace them by two
one-place maps. The reason is that according to the previous result, $s\sb {AG}\sp {(1)}(0,g)=0$ and $s\sb
{AG}\sp {(2)}(u,e)=e$. Thus we can define:
\begin{equation}\label{twomaps}
\pi \sb {AG }(a)=s\sb {AG}\sp {(1)}(a,e),\quad \varrho \sb {AG }(g)=s\sb {AG}\sp {(2)}(0,g)
\end{equation}
and obtain
\begin{equation}\label{prop(0,e)}
\pi \sb {AG }(0)=s\sb {AG}\sp {(1)}(0,e)=0,\quad \varrho \sb {AG }(e)=s\sb {AG}\sp {(2)}(0,e)=e.
\end{equation}
\begin{lem}\label{newFunction} If $\nu :(A,G,\cdot)\to (B,H,\bullet )$ is a homomorphism  and $\mu =\Phi (\nu )$,
then
$$\mu \sp {(1)} =\pi \sb {BH }\circ \nu \sp {(1)} \circ \pi \sb {AG }\sp {-1},$$
$$\mu \sp {(2)} =\varrho \sb {BH }\circ \nu \sp {(2)} \circ \varrho \sb {AG }\sp {-1}.$$
\end{lem}
\begin{proof}
Denote $\Phi (A,G,\cdot)$ by $\tilde{A}\tilde{G}$ and $\Phi (B,H,\bullet )$ by $\tilde{B} \tilde{H}$. Take
$\tilde{b}\in \tilde{B}$ and calculate:
\begin{multline}
\mu (\tilde{b},e)=s\sb {BH }\circ \nu  \circ s \sb {AG }\sp {-1}(\tilde{b},e)=s\sb {BH }\circ \nu (\pi \sb {AG
}\sp {-1}(\tilde{b}), e)=\\
=s\sb {BH }(\nu \sp {(1)}\circ \pi \sb {AG }\sp {-1}(\tilde{b}),e) =(\pi \sb {BH }\circ \nu \sp {(1)} \circ \pi
\sb {AG }\sp {-1}(\tilde{b}),e).
\end{multline}
Hence
$$\mu \sp {(1)}(\tilde{b}) =\pi \sb {BH }\circ \nu \sp {(1)} \circ \pi \sb {AG }\sp {-1}(\tilde{b}).$$
The proof of the second statement is similar.
\end{proof}
It is possible to strengthen the last result. We have an evident equality:
$$\a \sb {(w,g)}\circ \nu \sb {(0,x)}=\a  \sb {(0,g)} .$$
Allying $\Phi $, we obtain $\Phi (\a \sb {(w,g)})\circ \nu \sb {(0,x)}=\Phi (\a \sb {(0,g)})$. Hence $\a \sb
{s(w,g)}\circ \nu \sb {(0,x)}=\a \sb {s(0,g)}$. Under definition of the composition of homomorphisms, we obtain
$s\sp {(2)}(w,g)=(0,s\sp {(2)}(0,g))=(0,\varrho (g))$. Finally,
$$s(w,g)=(s\sp {(1)}(w,g),\varrho (g)).$$

 The obtained results gives us an opportunity to apply the second part of \thmref{potent-inner}
and conclude that $\Phi $ takes every free representation to an isomorphic one and assume in what follows that
$\Phi $ leaves fixed all objects of the category $(Rep-R )\sp 0$. Further, we can repeat all considerations in
the proof of \thmref{saveGenerators} and assume that $\Phi $ is such automorphism of $(Rep-R )\sp 0$ that
permutations $\pi \sb {AG }$ of $A$ and $\varrho \sb {AG }$ of $G$ leave fixed the basis of $(A,G,\cdot)$. This
gives us the following important fact
\begin{lem}\label{finalform}
$s\sb {AG} =\pi \sb {AG }\times \varrho \sb {AG }.$
\end{lem}
\begin{proof}Indeed,
$$s\sb {AG}(w,g)=\Phi (\a \sb {(w,g)})(1,x)=(\pi \sb {AG }\times \varrho \sb {AG })\circ \a \sb {(w,g)}
\circ (\pi \sb {W\sb 1 F\sb 1 }\times \varrho \sb {W\sb 1 F\sb 1 })\sp {-1} (1,x)=$$
$$=(\pi \sb {AG }\times \varrho \sb {AG })\circ \a \sb {(w,g)}(1,x)=(\pi \sb {AG }\times \varrho \sb {AG
})(w,g)=$$
$$=(\pi \sb {AG }(w),\varrho \sb {AG }(g)).$$

\end{proof}
\subsection{Derived binary operations} According to our method, we have to find the derived operations such that
the permutation $s$ is an isomorphism onto derived structure. Since the most arity of operations is 2, we
consider two-generated free representations. Let $(W\sb 2 ,F\sb 2)$ be a free two-generated representation,
where $ F\sb 2$ is the free group with two free generators $x\sb 1,x\sb 2$ and $W\sb 2 $ is the free $RF\sb 2
-$module with bases $Y=\{y\sb 1 , y\sb 2 \}$ : $W\sb 2 =y\sb 1 RF\sb 2 \oplus y\sb 2 RF\sb 2 $. Denote for
simple $\pi \sb {W\sb 2 F\sb 2 }$ by $\pi$ and $\varrho \sb {W\sb 2 F\sb 2 }$ by $\varrho $.
\begin{thm} The map $\varrho$ is the identity or the mirror map on $F\sb 2$, the last one means that
$\varrho$ assigns to every word in the free group $F\sb 2$ the same word written in the reverse order.
\end{thm}
\begin{proof} Since $\varrho (e)=e$, the derivative binary operation $x\sb 1\diamond x\sb 2 =\varrho (x\sb 1 x\sb 2)$
has the same unit $e$ and hence the same inverses. It implies that  $x\sb 1 \diamond x\sb 2=x\sb 1 x\sb 2$ or
$x\sb 1 \diamond x\sb 2=x\sb 2 x\sb 1$.
\end{proof}

Now it is turn to consider the derived additive operation on $W\sb 2$.  It is defined in the following way:
$$y\sb 1 \bot  y\sb 2 =\pi (y\sb 1 + y\sb 2).$$
Being an element of free module $W\sb 2$, $y\sb 1 \bot y\sb 2$ has a form $y\sb 1 P\sb 1 +y\sb 2 P\sb 2$, where
$P\sb 1 $ and $P\sb 1 $ are elements of group algebra $RF\sb 2$. Because of commutativity of the considered
operation, $P\sb 1 =P\sb 2 =P$ and $y\sb 1 \bot y\sb 2=(y\sb 1 +y\sb 2 ) P$. We know that $\pi (0)=0$. Therefore
$y\sb 1 =y\sb 1 P$ and hence $P=1$. Finally, we obtain that $y\sb 1 \bot y\sb 2=y\sb 1 +y\sb 2 $ and therefore
$\pi$ is an automorphism of the additive group of the module $W\sb 2$.

\begin{lem}\label{s-structure} The permutation $s$ of the set $W\sb 2 \times F\sb 2$ has the form: $s=\pi \times
\rho $, where $\pi $ is an automorphism of additive group of module $W\sb 2$ and $\rho $ is the identity map or
the mirror permutation of the group $F\sb 2$. For every $r\in R , w\in W\sb 2 $, it holds: $\pi (ru)=\varphi
(r)\pi (u)$, where $\varphi $ defined in \corref{MultRingMon} is automorphism of the ring $R$.
\end{lem}
\begin{proof}
We have to proof only the second statement. It is clear that $\a \sb {(ru, x\sb 1)}=\a \sb {(u, x\sb 1)}\circ
\nu \sb {(r,x)}$. Applying the automorphism $\Phi$, one obtain $\a \sb {(\pi (ru), x\sb 1)}=\a \sb {(\pi (u),
x\sb 1)}\circ \nu \sb {(\varphi(r),x)}$. Hence $\pi (ru)=\varphi(r)\pi (u)$. It was proved
(\corref{MultRingMon}), that $\varphi $ is automorphism of the multiplicative monoid of the ring $R$. Further,
$\pi ((r+t)y\sb 1)=\varphi (r+t)y\sb 1$ and $\pi ((r+t)y\sb 1)=\pi (r y\sb 1 +t y\sb 1)=\varphi (r) y\sb 1
+\varphi(t) y\sb 1 $. Hence $\varphi (r+t)y\sb 1= (\varphi (r) +\varphi(t)) y\sb 1 $ and finally $\varphi
(r+t)=\varphi (r)+\varphi (t)$.
\end{proof}
Let $\varphi$ be an automorphism of the ring $R$. It can be extended up to so called {\it twisted} automorphism
of every free $RF-$module $W$ for a group $F$ as follows: $\varphi \sb {WF} (y\sum \sb g r\sb g g)=y\sum \sb g
\varphi (r\sb g) g$.
\begin{defn} Consider the function which assigns to every free representation $(W,F,\cdot )$ the pair $(\varphi \sb
{WF},1\sb F)$ of permutations, where $\varphi \sb {WF}$ is the twisted automorphism of $W$ defined above and
$1\sb F$ is the identity on $F$. This function determines an automorphism of the category $(Rep-R )\sp 0$. We
call this automorphism a {\it standard twisted} automorphism and denote it by $\hat{\phi}$.
\end{defn}

\subsection{Derived action.}
Now consider the representation structure, that is the action of the group $F\sb 2 $ on the module $W\sb 2$. We
assume that the ring $R$ without zero devisors. The action is determined by the term $y\sb 1 \cdot x\sb 1$. The
derived structure is determined by the term $y\sb 1 \bullet x\sb 1 =\pi (y\sb 1 \cdot x\sb 1)$ and the
permutation $s=\pi \times \rho $ is a isomorphisms between this two structures. Let $\nu \sb {(x, e)}$ be the
endomorphism of $(W\sb 1, F\sb 1$) defined as usual by $\nu \sb {(x, e)} (1,x)=(x,e)$. We have
$$\a \sb {(y\sb 1,x\sb 1)}\circ \nu \sb {(x, e)}=\a \sb {(y\sb 1 \cdot x\sb 1 ,e)}.$$
Apply our automorphism $\Phi $ and obtain
$$\a \sb {(y\sb 1,x\sb 1)}\circ \nu \sb {(w, e)}=\a \sb {(y\sb 1 \bullet x\sb 1 ,e)},$$
where $w$ is an element of the group algebra $RF\sb 1$, that is $w=\sum r\sb i x\sp i ,\; i\in \mathbb{Z}$.
Hence
$$y\sb 1 \bullet x\sb 1 =y\sb 1 \cdot \sum r\sb i x\sb 1\sp i . $$
and fir every $u\in W\sb 2 $ and every $g\in F\sb 2$ we have
$$u \bullet g =u \cdot \sum r\sb i g\sp i . $$
Since the derived structure is isomorphic to the source one, $y\sb 1 \bullet e =y\sb 1$, that gives $\sum r\sb i
=1$. Write the associativity low of the derived action and obtain
$$(y\sb 1 \bullet x\sb 1 )\bullet x\sb 2=y\sb 1 \bullet (x\sb 1 x\sb 2),$$
if $\rho $ is the identity map and
$$(y\sb 1 \bullet x\sb 1 )\bullet x\sb 2=y\sb 1 \bullet (x\sb 2 x\sb 1),$$
if $\rho $ is the mirror map.  Consider the first case and obtain the identity in the variety of all group
representation:
$$(\sum r\sb i x\sb 1 \sp i )(\sum r\sb i x\sb 2 \sp i) = \sum r\sb i (x\sb 1 x\sb 2) \sp i .$$
the both sides of this equality can have only two pairs common words: one of degree 0 and one of degree 2, that
is, $x\sb 1 x\sb 2$.  It implies that $r\sb i =0 $ for all $i\not = 1$ and $r\sb 1 =1$. That is $y\sb 1 \bullet
x\sb 1 =y\sb 1 \cdot x\sb 1$. It is easy to see that in the case $\rho $ is the mirror map we have $y\sb 1
\bullet x\sb 1 =y\sb 1 \cdot x\sb 1 \sp {-1}$.
\begin{cor}\label{pi_action} For every $u\in W\sb 2$ and $g\in F\sb 2$, $\pi (u\cdot g)=\pi (u)\cdot g$ if $\rho $ is the identity,
and $\pi (u\cdot g)=\pi (u) \cdot \bar{g}\sp {-1}$, if $\rho$ is the mirror map.
\end{cor}
\begin{proof} Under definition, $\pi (u\cdot g)=\pi (u)\bullet \rho (g)$. Thus for the case $\rho $ is the
identity, we obtain the first statement. If $\rho $ is the mirror map $\pi (u)\bullet \rho (g)=\pi (u)\cdot
\bar{g}\sp {-1}$.

\end{proof}

This fact leads to a definitions of a new kind of quasi-homomorphisms. For every free group $F$, the map
$g\mapsto \bar{g} $ can be extended up to a permutation of the group algebra $RF$ as follows: $\overline{rg} =
r\bar{g}$ for every $r\in R$ and $g\in F$. In the same way, the map $g \mapsto g\sp {-1}$ can be extended up to
$RF$. Clearly, these two maps can be extended up to every free $RF-$module $W$ as follows: $(yP) \sp {-1} =yP\sp
{-1}$ for every free generator $y$ and $P\in RF$. The same for "bar".
\begin{defn}\label{newmirror}
Let $(W, F, \cdot )$ is a free representation. The pair $\delta =(\delta \sp {(1)},\delta \sp {(2)}) $, where
$\delta \sp {(1)}(w)=\bar {w}\sp {-1}$ for $w\in W$ and $\delta \sp {(2)}(g)=\bar{g}$ for $g\in F$, is called a
mirror automorphism of $(W, F, \cdot) $. Corresponding to every free representation $(W, F, \cdot )$ the mirror
automorphism $\delta \sb {WF}$ we obtain so called {\it mirror} automorphism $\Delta $ of the category $(Rep-R
)\sp 0$.
\end{defn}
\subsection{Final.} Let $R$ be a ring without zero devisors. Now we are able to describe automorphisms of the
category $(Rep-R )\sp 0$ taking the regular representation $(W\sb 1, F\sb 1)$  to an isomorphic one in the
similar way as in the previous section.

\begin{thm}\label{QRep-R} Let $\Q $ be an extension of the category $(Rep-R )\sp 0$ obtained by adding twisted
and mirror automorphisms. Then every automorphism of $(Rep-R )\sp 0$ is $\Q -$inner.
\end{thm}
It turns out that this description can be simplified. To show it consider a function ${\bf c}:(W, F, \cdot
)\mapsto c\sb {WF}$ which assigns to every free group representation $W, F, \cdot )$ the permutation $ c\sb
{WF}$ as follows:
$$c\sb {WF}(w,g)=(w,g\sp {-1})$$
for every $w\in W $ and $g\in F$.
\begin{lem}The function $\bf c $ is central for the category $(Rep-R )\sp 0$ and every map $c\sb {WF}$ is an
isomorphism of $(W, F, \cdot )$ onto the derivative structure $(W, F\sp *, \bullet )$, where $w\bullet g=w\cdot
g\sp {-1}$ and $F\sp *$ is the dual to $F$ group.
\end{lem}
\begin{proof} Let $\mu :(A, G, \cdot )\to (B, H, \ast )$ be a homomorphism of in the category $(Rep-R )\sp 0$.
Calculate for every $a\in A $ and $g\in G$:
$$\mu \sb {AG} \circ c\sb {AG}(a,g)=\mu \sb {AG}(a,g\sp {-1})=(\mu \sb {AG}\sp {(1)}(a),(\mu \sb {AG}\sp
{(2)}(g))\sp {-1})=$$
$$=c\sb {BH}(\mu \sb {AG} \sp {(1)}(a),\mu \sb {AG} \sp {(2)}(g)))=c\sb {AG}\circ \mu \sb {BH}(a,g).$$
Thus the function $\bf c $ is central for the category $(Rep-R )\sp 0$.

Further, since the map $c\sb {WF}\sp {(1)}$ is the identity and the map $c\sb {WF}\sp {(2)}: g\mapsto g\sp {-1}$
is an isomorphism of the group $F$ onto the dual group $F\sp *$, we have only to prove that $c\sb {WF}$ is a
homomorphism with respect to group actions. For every $w\in W $ and $g\in F$ we have:
$$ c\sb {WF}\sp {(1)}(w\cdot g)=w\cdot g =w \bullet g\sp {-1}=c\sb {WF}\sp {(1)}(w)\bullet c\sb {WF}\sp
{(2)}(g).$$ Thus $c\sb {WF}$ is an isomorphism of $(W, F, \cdot )$ onto the derivative structure $(W, F\sp *,
\bullet )$.
\end{proof}
Applying \lemref{InnerAndCentral}, we obtain that mirror automorphisms of the category $(Rep-R )\sp 0$ are in
fact inner. Thus the final description is the following

\begin{thm}\label{Rep-R}. Let an automorphism $\Phi$ of $(Rep-R )\sp 0$  take the regular representation
$(W\sb 1, F\sb 1)$ to an isomorphic one. Then  $\Phi$ is a semi-inner automorphism, that is, it can be
represented in the following form $ \Phi = \hat{\varphi}\circ \Psi ,$ where $\Psi $ is inner and
$\hat{\varphi}$ is the standard twisted $\varphi -$automorphism for some automorphism $\varphi$ of the ring $R$.
\end{thm}
\begin{proof} $\Phi $ is a composition of an inner automorphism $\Psi$ and an automorphism $\Phi \sb 1$
that leaves fixed all objects and its bases. According to \lemref{s-structure}, $\Phi \sb 1$ determines an
automorphism $varphi$ of the ring $R$ and hence the standard twisted automorphism $\hat{\varphi} $ of our
category. Then the automorphism $\Gamma =\Phi \sb 1 \circ \hat{\varphi} \sp {-1}$ determines the identity
automorphism of the ring $R$ and according to \corref{pi_action} it is the identity automorphism.
\end{proof}

The property of the category $(Rep-R )\sp 0$ that its automorphisms take the regular representation to an
isomorphic one depends on the ring $R$. Since this condition is satisfied if $R$ is an infinite field $\K$, all
automorphisms of the category $(Rep-\K )\sp 0$ are semi-inner.

Institute of Mathematics, Hebrew University, 91803 Jerusalem, Israel

{\it E-mail address}: borisov@math.huji.ac.il

Department of Mathematics, Bar-Ilan University, 52900 Ramat Gan, Israel

{\it E-mail address}: zhitomg@macs.biu.ac.il

\end{document}